\documentclass[MSNbibl,number,citesort,dvips]{arxbj}
\usepackage{upgreek}

\aid{0}
\volume{18}
\issue{3}
\pubyear{2012}
\firstpage{857}
\lastpage{882}
\doi{10.3150/11-BEJ366}

\makeatletter
\newcommand{\reals}{{\mathbb R}}
\newcommand{\bfB}{{\mathbf B}}
\newcommand{\bfJ}{{\mathbf J}}
\newcommand{\bfT}{{\mathbf T}}
\newcommand{\bfu}{{\mathbf u}}
\newcommand{\bfx}{{\mathbf x}}
\newcommand{\intTh}{{\operatorname{int}(\Theta)}}

\newcommand{\nmk}{{\biggl(\frac{n}{n-m}\biggr)^{k/2}}}
\newcommand{\nmkalt}{{\biggl(1+\frac{mk}{2n}\biggr)}}
\newtheorem{theorem}{Theorem}
\newtheorem{lemma}{Lemma}

\renewcommand{\epsilon}{\varepsilon}
\newcommand{\cal}{\mathcal}
\newproclaim{condition}{Condition}
\newremark{remark}{Remark}
\newcommand{\fraca}[2]{{#1}/{#2}}
\makeatother

\begin{document}
\begin{frontmatter}

\title{Conditional limit laws for goodness-of-fit tests}
\runtitle{Conditional limit laws}

\begin{aug}
\author{\fnms{Richard A.} \snm{Lockhart}\corref{}\ead[label=e1]{lockhart@sfu.ca}\ead[label=u1,url]{www.stat.sfu.ca/\textasciitilde lockhart}}
\runauthor{R.A. Lockhart}
\address{Department of Statistics and Actuarial Science,
Simon Fraser University,
Burnaby, BC V5A 1S6,
Canada. \printead{e1,u1}}
\end{aug}

\received{\smonth{5} \syear{2010}}
\revised{\smonth{1} \syear{2011}}

%
\begin{abstract}
We study the conditional distribution of goodness of fit statistics of
the Cram\'{e}r--von Mises type given the complete sufficient statistics
in testing for exponential family models. We show that this distribution
is close, in large samples, to that given by parametric bootstrapping,
namely, the unconditional distribution of the statistic under the value
of the parameter given by the maximum likelihood estimate.
As part of the proof,
we give uniform Edgeworth expansions of Rao--Blackwell estimates in these
models.
\end{abstract}

%
\begin{keyword}
\kwd{empirical distribution function}
\kwd{goodness-of-fit}
\kwd{local central limit theorem}
\kwd{parametric bootstrap}
\kwd{Rao--Blackwell}
\end{keyword}

\end{frontmatter}
%

\section{Introduction}\label{sec1}

In this paper, we compare conditional and unconditional goodness-of-fit tests
and give conditions under which the two give essentially identical
results in large samples. Our results apply in testing fit
for exponential family models for independent and identically
distributed (i.i.d.)
data, $X_1,\ldots,X_n$. Our interest is to test the null hypothesis
that the distribution of the individual $X_i$ belongs
to a natural exponential family with density, relative to some
$\sigma$-finite measure, $\mu(\mathrm{d}x)$, on some sample space $\Omega$,
of the form
%
\begin{equation}\label{expdensity}
f(x;\theta) \equiv c(x) \exp\{\theta^\prime T(x)-\kappa(\theta)\}
\end{equation}
with natural parameter space $\Theta\subset\reals^k$; we assume that
$\Theta$
has non-empty interior which we denote $\intTh$. In (\ref
{expdensity}), $T$
takes values in $\reals^k$ and superscript $\prime$ denotes transposition.
A complete and
sufficient statistic for the parameter $\theta$ is then
\[
\bfT_n\equiv\bfT_n(X_1,\ldots,X_n) = \sum_{i=1}^n T(X_i).
\]

To apply classical hypothesis testing ideas, we regard this model as a null
hypothesis. We consider the omnibus alternative hypothesis that the
sample is drawn from a distribution which is not in the parametric model.
One common approach to this hypothesis testing problem is to define some
statistic $S(X_1,\ldots,X_n;\theta)$ which measures in some way departure
of the sample from what is expected if $\theta$ is the true value. Since
$\theta$ is unknown, it is replaced in this measure by $\hat\theta
_n$, the
maximum likelihood estimate of the parameter vector, leading to the
statistic $S_n\equiv S(X_1,\ldots,X_n,\hat\theta_n)$.

Common examples include empirical distribution function
statistics such as Cram\'er--von Mises, Kolmogorov--Smirnov, Anderson--Darling
and many chi-squared statistics. The usual situation is
that the test statistic has a distribution which depends, even in large
samples, on the unknown parameter value (exceptions arise in the normal and
other families which have only location and/or scale
parameters). Thus, to implement the tests in
practice it is necessary to specify how to compute critical points for the
tests or how to compute appropriate $P$-values corresponding to the test
statistics.
A~method long in use is to derive large sample theory for the statistic
$S_n$, establishing the convergence in distribution of $S_n$ to some
limiting distribution which depends on the true value of~$\theta$. If
$C_\alpha(\theta)$ is the upper $\alpha$ critical point of this limiting
distribution and $C_\alpha$ depends continuously on $\theta$, then
the test
which rejects if $S_n > C_\alpha(\hat\theta_n)$ has asymptotic level
$\alpha$.
See Lockhart and Stephens \cite{LSvM} for a discussion of this method in testing fit
for the von Mises distribution for directional data; this testing problem
is discussed below in more detail.

A more modern method which achieves the same asymptotic behaviour is
the parametric bootstrap. Let $H_n(\cdot;\theta)$ denote the cumulative
distribution function of $S_n$ when the true
parameter value is $\theta$.
Then
\[
P_b = 1-H_n(S_n;\hat\theta_n)
\]
is the parametric bootstrap $P$-value.
This $P$-value is usually computed approximately by generating some number,
$B$, of bootstrap samples drawn from the density $f(\cdot,\hat\theta_n)$,
computing the statistic $S_n$ for each of these $B$
samples and then counting the fraction of these bootstrap statistic values
which exceed the value of $S_n$ for the data set at hand.

These two methods for goodness-of-fit testing both depend on asymptotic
theory to justify their performance. They do not have, except in the
location-scale situation mentioned, exact level $\alpha$ and thus no
exact finite sample optimality properties. Conditional tests, which we
discuss next, offer at least the potential for such optimality. (See
Remark~\ref{rem9} in the \hyperref[discussion]{Discussion} section for some comments.)

One standard approach (discussed in detail in   \cite{LR}) to optimality
theory is to search for powerful unbiased level $\alpha$ tests: tests whose
power never falls below $\alpha$ on the alternative. Such tests will
generally have Neyman structure; that is, their level will
be $\alpha$ everywhere on the boundary of the null hypothesis.
For the omnibus alternative, this boundary is generally the entire model.

Now suppose $\bfT_n$ is a complete sufficient statistic for
this model. Then the requirement that the level of the test be $\alpha$
everywhere in the parametric model and completeness guarantee that the
test must have conditional level $\alpha$. That is, an
unbiased level $\alpha$ test must have the property that the conditional
probability of
rejection given $\bfT_n$ is identically $\alpha$. This is precisely the
argument used in
Lehmann and Romano \cite{LR} to show that Student's $t$ test is uniformly most powerful unbiased.

By a conditional test, then, we mean a test whose level, given the
sufficient statistic $\bfT_n$, is identically $\alpha$. Two recent
papers on
goodness-of-fit, Lockhart, O'Reilly and Stephens~\cite{LOS1,LOS2},
have compared such conditional tests with parametric bootstrap
tests.
They implemented their conditional tests as follows. For a test
statistic $S_n$, let
$G_n(\cdot|\cdot)$ denote the conditional distribution function, when the
true distribution of the data comes from the exponential family,
of $S_n$ given $\bfT_n$.
This function $G_n$ does not depend on $\theta$.
If this conditional distribution function is continuous then
\[
P_c=1-G_n(S_n|\bfT_n)
\]
has a uniform distribution under $H_o$; it is therefore an exact $P$-value.
These $P$-values are often computed by Monte Carlo or Markov Chain Monte
Carlo; see Lockhart, O'Reilly and Stephens \cite{LOS1,LOS2}
for examples and references.

In Lockhart, O'Reilly and Stephens \cite{LOS2}, for instance, the authors considered an i.i.d. sample
from the von Mises distribution. Observations $X_i$ are points on
the unit circle; see Section \ref{vonmisessec} below
for details of the density. The complete sufficient statistic is
$\bfT_n=\sum X_i$ and
the authors
use Watson's $U^2$ statistic for $S_n$. They use
Markov Chain Monte
Carlo methods to generated a sequence of samples from the conditional
distribution of
$X_1,\ldots,X_n$ given $\bfT_n$;
all the generated samples have the same value of $\bfT_n$.
The authors evaluate $P_c$ by computing $U^2$ for
each data set and estimating $P_c$ by the fraction of samples giving
larger values of $U^2$ than the original data
sample.

These authors also compute the parametric bootstrap value, $P_b$,
for the same statistic by generating i.i.d. samples from
the von Mises distribution using, for the parameter
value, the estimate
of the parameter derived from the original data. Of course the values
of the sufficient statistic $\bfT_n$ vary from one bootstrap sample to
another. Again $U^2$ is computed for each
bootstrap sample and a $P$ value is computed as the fraction of
bootstrap $U^2$ values which are larger
than the observed value of $U^2$.

Very high correlations between the $P$-values computed using these two
methods were observed in Lockhart, O'Reilly and Stephens \cite{LOS2}. For example, they considered a
test that a sample of size 34 comes from a von Mises distribution.
Using Watson's $U^2$ and generating samples from the null
hypothesis they observed a correlation of 0.997
between the two $P$-values. For a sample of 55 observations, the
correlation observed was 0.9997.

Here we show that for statistics $S_n$ of the
Cram\'{e}r--von~Mises type these two methods must give similar
$P$-values because, when the null hypothesis is true,
\[
\sup_s\{|G_n(s|\bfT_n)-H_n(s;\hat\theta_n)|\} \to0
\]
in probability, at least when the model being tested is an exponential family.
In fact, the convergence is almost sure for samples from any distribution
for which $\hat\theta_n/n$ converges almost surely to an interior point
of the parameter space. For statistics $S_n$ which are sums of the form
$\sum_i u_n(X_i,\hat\theta)$ this result is established by Holst \cite{holst}.
Our results extend his to statistics which we now describe.

When $\Omega$ is the real line, many goodness-of-fit tests are based on
statistics which are functionals of the estimated empirical process
\[
W_n(s) = \sqrt{n} \{F_n(x) -F(x,\hat\theta_n) \},
\]
where we now use $F(x,\theta)$ for the cumulative distribution function,
$s$ is related to $x$ by
$ s=F(x,\hat\theta_n)$
and
$F_n$ is the usual empirical distribution function:
\[
F_n(x) = n^{-1}\sum_{i=1}^n 1(X_i \le x).
\]
Common choices for statistics include:
\begin{itemize}
\item Cram\'{e}r--von Mises type:
%
\begin{equation}\label{cvmdef}
S_n = \int_{0}^1 \psi^2(s) W_n^2(s)\,\mathrm{d}s;
\end{equation}
\item Watson type:
\[
\int_{0}^1  \biggl\{W_n(s) - \int_{0}^1\psi(u)W_n(u)\,\mathrm{d}u  \biggr\}^2
\psi^2(s)\,\mathrm{d}s;
\]

\item Kolmogorov--Smirnov type:
\[
\sup_{0 < s < 1}  |\psi(s)W_n(s) |.
\]
\end{itemize}
In each case, $\psi$ is some weight function defined on $(0,1)$.

The large sample analysis of the unconditional distribution of
such statistics comes from the well known weak convergence, in $D[0,1]$,
of the process $W_n$ to a Gaussian process, $W$, which we now describe.
Let ${\cal I}(\theta)$ be the Fisher information matrix and define the
column vector
\[
\xi(s,\theta)= \frac{\partial F(x;\theta)}{\partial\theta},
\]
where $x$ is defined as a function of $s$ by
$
F(x,\theta)\equiv s.
$
Then the limit process $W$ has mean~0 and covariance function
\[
\rho_\theta(s,t) =
\min\{s,t\}-st
-\xi(s,\theta)^\prime \{{\cal I}(\theta) \}^{-1} \xi
(t,\theta).
\]
The statistics indicated above are all continuous functionals of
$W_n$ (under mild conditions on the weight functions involved) and
as such converge in distribution to the same functional applied to
the limit process $W$. See Stephens \cite{stephensEDFch4}
for a detailed discussion
of the resulting tests and Shorack and Wellner \cite{shorackwellner}
for mathematical details.

The weak convergence result can be proved in two steps:
prove convergence in distribution of the
finite dimensional distributions of $W_n$ and then prove
tightness of the sequence of processes in $D[0,1]$.
We believe a similar result holds, in exponential families,
conditional on the sufficient statistic. Results in Holst \cite{holst}
can be used to establish convergence of the conditional finite
dimensional distributions but we are unable to extend
the calculations to prove conditional tightness. Instead we use
Holst's
results and a~truncation argument to
deal directly with statistics of the Cram\'{e}r--von Mises or Watson types.
Without tightness we cannot handle statistics of the Kolmogorov--Smirnov
type.

Our truncation argument uses an accurate approximation to the conditional
expectation, given $\bfT_n$, of the statistic in question. This approximation
is based on an expansion of the difference between a Rao--Blackwell
estimate and
the corresponding maximum likelihood estimate. Our results here extend the
work of Portnoy \cite{portnoy}.

Section~\ref{sec2} gives precise statements of our results for the case
of Cram\'er--von Mises statistics. Section~\ref{sec3} gives the expansion of the
Rao--Blackwell estimate. Section~\ref{sec4} applies the calculations to two examples
showing how to verify the main condition, Condition~\ref{condD} below,
and illustrating the expansions of Section~\ref{sec3}.
Section~\ref{sec5} provides some discussion and indicates the
extension to Watson's statistic and other statistics which are quadratic
functionals of the empirical distribution. In that section,
we consider power and discuss various rephrasings of our main result.
Details of some proofs are in Section~\ref{proofs}.

\section{Main results}\label{sec2}

\subsection{Absolutely continuous distributions}\label{sec2.1}

We seek to test the hypothesis that the distribution
of each $X_i$ belongs
to a natural exponential family with density, relative to
some $\sigma$-finite measure
$\mu(\mathrm{d}x)$ on $\Omega$,
of the form (\ref{expdensity}) and complete sufficient
statistic $\bfT_n$ as described in the \hyperref[sec1]{Introduction}.
We will need a number of well known facts about exponential families
which we
gather here in the form of a lemma.

\begin{lemma}\label{exponentialfamilyprops}
The random vector $\bfT_n$ has
moment generating function
\[
{\rm E}_\theta [\exp\{\phi^\prime\bfT_n\} ] =
\exp [n\{\kappa(\phi+\theta)-\kappa(\theta)\} ]
\]
which is finite whenever $\theta+\phi\in\Theta_0$, and
cumulants $n\kappa_{i_1,\ldots,i_r}$ where
\[
\kappa_{i_1,\ldots,i_r}
= \frac{\partial^r \kappa(\theta)}{\partial\theta_{i_1}\,\cdots\, \partial\theta_{i_r}}.
\]
In particular, the mean of $\bfT_n$ is
\[
{\rm E}_\theta(\bfT_n) = n\mu(\theta)\equiv n\nabla\kappa(\theta),
\]
where $\nabla$ is the gradient operator. The covariance matrix is
\[
\operatorname{Var}_\theta(\bfT_n) \equiv nV(\theta) = n\nabla^2\kappa
(\theta),
\]
where $\nabla^2$ denotes the Hessian operator. Thus, $V(\theta)$ has
entries
\[
V_{ij}(\theta) =
\frac{\partial^2\kappa(\theta)}{\partial\theta_i\,\partial\theta_j}.
\]
Moreover, all moments and cumulants of $\bfT_n$ depend smoothly on
$\theta$ on
the interior of $\Theta$.
\end{lemma}

Our results apply to exponential families where $\bfT_n$
has a density relative to Lebesgue measure. We assume the following condition.

\renewcommand{\thecondition}{D}
\begin{condition}\label{condD} For every compact subset $\Gamma$ of
$\intTh$, there is an integer $r$ such that the characteristic function
\[
\eta_\theta(\phi) \equiv{\rm E}_\theta \{\exp(\mathrm{i}\phi^\prime
\bfT_r) \} =
\exp [r\{\kappa(\theta+\mathrm{i}\phi)-\kappa(\theta)\} ]
\]
is integrable for all $\theta\in\Gamma$ and
\[
\sup_{\theta\in\Gamma}\int_{\reals^k} |\eta_\theta(\phi)|\,\mathrm{d}\phi < \infty.
\]
\end{condition}

Condition~\ref{condD} has two consequences we need. First, it means the
matrix
$
\operatorname{Var}_\theta(\bfT_1) = \nabla^2\kappa(\theta)
$
is positive definite for each $\theta\in\intTh$.
This
implies the map $\theta\mapsto\mu(\theta)=\nabla\kappa(\theta)$
is an
open bijective mapping of $\intTh$ to $\mu(\intTh)$. A second
consequence is that $\bfT_n$ has bounded continuous density
for each $\theta\in\Gamma$ and $n\ge r$. In the
examples it will be useful to know the converse is also true. The
following lemma is essentially Theorem~19.1 in
Bhattacharya and Ranga~Rao \cite{bhatta}, page~180;
see also Lemma~\ref{HolstLemma} in Section~\ref{proofs} below.

\begin{lemma} \label{integrability} Condition~\textup{\ref{condD}} is equivalent
to Condition~\textup{\ref{condD*}}.
\end{lemma}

\renewcommand{\thecondition}{D$^*$}
\begin{condition}\label{condD*} For every compact subset $\Gamma$
of $\intTh$ there is an integer $r$ such that $\bfT_r$ has
continuous (Lebesgue) density $f_r(t;\theta)$ for each $\theta\in
\Gamma$ and
\[
\sup_{\theta\in\Gamma}\sup_{t\in\reals^k} f_r(t,\theta) <
\infty.
\]
\end{condition}

As in the \hyperref[sec1]{Introduction}, we let
$
G_n(\cdot| t)
$
denote the conditional cumulative distribution function of $S_n$
given $\bfT_n=t$. Also let
$
H_n(\cdot; \theta)
$
denote the unconditional cumulative distribution function
of $S_n$ when $\theta$ is the value of the parameter.

We will show that for statistics which are sums as in (\ref{sumdef}) below
or of the Cram\'er--von~Mises type these two cumulative distributions are
uniformly close provided that $t$ and $\theta$ are related properly, that is,
$t=n\mu=n\nabla\kappa(\theta)$.

Our results use a minor modification of Corollary 3.6 of Holst \cite{holst}
which establishes this uniform closeness for statistics which are sums over
the data as described below. We use the following
notation. By ${\cal L}(S_n;\theta)$
we mean the unconditional
distribution of $S_n$ under the model with true parameter $\theta$. By
${\cal L}(S_n|\bfT_n=t)$
we mean the conditional
distribution of $S_n$ given $\bfT_n=t$. We use the symbol $\Rightarrow
$ to
denote convergence in distribution (weak convergence) and ${\cal L}(W)$ and
similar notation for limiting distributions.
Our version of Holst's results is:

\begin{lemma} \label{Holst} Assume Condition~\textup{\ref{condD}}.
Suppose that $u_n(\cdot;\cdot)$ is a sequence of measurable functions
mapping $\Omega\times\Theta$ to $\reals^m$.
Let
%
\begin{equation}\label{sumdef}
S_n(\theta) = n^{-1/2} \sum_{i=1}^n  [ u_n(X_i,\theta) - {\rm
E}_{\theta} \{u_n(X_i,\theta) \} ].
\end{equation}
%
Assume that for any deterministic sequence $\theta_n$ of parameter values
converging to some $\theta\in\intTh$ the joint law
\[
{\cal L}_{\theta_n} \bigl(S_n(\theta_n),
n^{-1/2}  \{\bfT_n - n\mu(\theta_n) \} \bigr)
\]
converges to
multivariate normal with mean 0 and variance--covariance matrix
of the form
\[
 \left[
\matrix{ A(\theta) & B(\theta) \cr
B^\prime(\theta) &
V(\theta)
}
 \right]
\]
which may depend on $\theta$ but not on the specific sequence
$\theta_n$.
Then with $S_n$ denoting $S_n(\hat\theta_n)$ we
have for every such sequence $\theta_n$
\[
{\cal L}(S_n|T_n=t_n) \equiv G_n(\cdot| n\mu(\theta_n) )\Rightarrow \operatorname{MVN}\bigl(0,
A(\theta)-B(\theta)V^{-1}(\theta)B^\prime(\theta)\bigr),
\]
where $t_n =n\mu(\theta_n)$.
Moreover, for every compact subset $\Gamma$ of
$\intTh$ we have
\[
\lim_{n\to\infty}\sup_{-\infty<x<\infty}
\sup_{\theta\in\Gamma} |G_n(x|n\mu)-H_n(x|\theta)| = 0.
\]
\end{lemma}

The condition involving the sequence $\theta_n$
amounts to requiring that the central limit theorem
apply uniformly on compact subsets of $\Theta$.
Our main result extends the last conclusion of the lemma
to statistics of the Cram\'er--von~Mises type for the case where
$\Omega$ is the real line; see Remark~\ref{rem8} in Section~\ref{discussion} for
discussion of more general sample spaces.

\begin{theorem} \label{maintheorem}
Suppose $S_n$ is as defined in (\ref{cvmdef}). Suppose the
weight $\psi$ is
continuous on $[0,1]$.
Assume Condition~\textup{\ref{condD}}.
Then for every compact subset $\Gamma$ of
$\intTh$ we have
%
\begin{equation}\label{mainresult}
\lim_{n\to\infty}\sup_{-\infty<x<\infty}
\sup_{\theta\in\Gamma} |G_n(x|n\mu)-H_n(x|\theta)| = 0.
\end{equation}
\end{theorem}

The theorem asserts that two distribution functions, one conditional,
the other unconditional, are close together
everywhere and simultaneously for all $\theta$ belonging to some compact
set. In the Introduction, we described our results in terms of $P$-values;
we now recast the theorem in those terms. The conditional $P$ value,
now denoted $P_{c,n}$, is
\[
P_{c,n} \equiv1-G_n(S_n|\bfT_n).
\]
The unconditional $P$ value, $P_{u,n}$, is
\[
P_{u,n} \equiv1-H_n(S_n;\hat\theta_n).
\]
We then have the following result which also clarifies the sampling
properties of the distributions $G_n$ and $H_n$ evaluated at sample
estimates.

\begin{theorem} \label{corr}
Assume the conditions of Theorem \ref{maintheorem}.
\begin{enumerate}[(b)]
\item[(a)] If $X_1,X_2,\ldots$ is an i.i.d. sequence generated from
the model
with true parameter value $\theta\in\intTh$ (i.e.,
if the null hypothesis is true and the true parameter value is not on the
boundary of the parameter space), then
\[
\lim_{n\to\infty}\sup_{-\infty<x<\infty}
\sup_{\theta\in\Gamma} |G_n(x|\bfT_n)-H_n(x|\hat\theta_n)| = 0
 \qquad \mbox{almost surely}
\]
and
\[
P_{c,n}-P_{u,n} \to0 \qquad \mbox{almost surely}.
\]
\item[(b)] Suppose $X_1,X_2,\ldots$ is an i.i.d. sequence generated from
some fixed alternative distribution. Suppose that for this
alternative ${\rm E}(\bfT_1) = \mu_a$ exists and is in the open
set $\mu(\intTh)$, that is, the image of the interior
of $\Theta$ under the map $\theta\mapsto\mu$. Then both conclusions of
part (a) still hold.
In particular, if one test is consistent against the
alternative then so is the other.
\end{enumerate}
\end{theorem}

Details of proofs are in Section \ref{proofs} but here we outline
the strategy of proof for our Theorem~\ref{maintheorem}.
Fix a complete orthonormal system of functions $g_j$ defined on
$[0,1]$; for definiteness we take $g_j(s)=\sqrt{2} \sin(\uppi j s)$.
Define
\[
U_{n,j} = \int_0^1 \psi(s) W_n(s) g_j(s)\,\mathrm{d}s.
\]
Then by Parseval's identity
\[
S_n = \sum_{j=1}^\infty U_{n,j}^2.
\]
The proof then has the following steps:
\begin{enumerate}[9.]
\item\label{stdweakconv} The sequence of distribution functions
$H_n(\cdot|\theta)$ converges weakly to a limiting distribution function
$H_\infty(\cdot| \theta)$; the convergence is uniform on compact
subsets of
$\Theta$. The distribution in question is the law of
\[
\int_0^1 \psi^2(s) W^2(s)\,\mathrm{d}s =\sum_{j=1}^\infty
U_{\infty,j}^2,
\]
where we define
\[
U_{\infty,j}=\int_0^1 \psi(s)W(s) g_j(s)\,\mathrm{d}s.
\]
This 
reduces the problem to proving that the sequence
$G_n(\cdot|n\mu)$ converges uniformly on compact subsets of $\intTh$
to $H_\infty(\cdot|\theta)$ where $\mu=\nabla\kappa(\theta)$.

\item
Uniform convergence is established by
considering an arbitrary sequence $\theta_n$ of parameter values
converging to some $\theta\in\intTh$ and showing that, with
$\mu_n = \nabla\kappa(\theta_n)$,
%
\begin{equation}\label{weakconvergenceb}
\lim_{n\to\infty}\sup_{-\infty<x<\infty}
|G_n(x | n\mu_n)-H_\infty(x;\theta)| = 0.
\end{equation}
\item
\label{uncondfdds}
Apply standard weak convergence ideas to see that for each $K$
fixed 
%
\[
{\cal L} ((U_{n,1},\ldots,U_{n,K});\theta_n )
\Rightarrow
{\cal L} ((U_{\infty,1},\ldots,U_{\infty,K}) ).
\]

\item\label{condfdds}
Use Holst's results to prove that
\[
{\cal L} \bigl((U_{n,1},\ldots,U_{n,K})|\bfT_n = n\mu_n
\bigr)\Rightarrow
{\cal L} ((U_{\infty,1},\ldots,U_{\infty,K}) );
\]
this is the same joint limit law as in the previous step.

\item Prove the sequence ${\cal L}(S_n|\bfT_n=n\mu_n)$
of conditional distributions of $S_n$ is tight.

\item\label{uncondtrunc}
Prove that there is a sequence $K_n$ tending to infinity sufficiently
slowly that
\[
{\cal L} \Biggl(\sum_{j=1}^{K_n} U_{n,j}^2 ;\theta_n \Biggr)
\Rightarrow
{\cal L} \biggl(\int_0^1 \psi^2(s) W^2(s)\,\mathrm{d}s \biggr).
\]
\item\label{condtrunc}
Prove the corresponding conditional result given $\bfT_n=n\mu_n$.
\item\label{tails}
Prove that for any sequence $K_n$ tending to infinity
\[
\sum_{j=K_n}^\infty U_{n,j}^2
\]
converges to 0 in probability given $\bfT_n=n\mu_n$.

\item Apply Slutsky's theorem to \ref{uncondtrunc},
\ref{condtrunc} and \ref{tails} and use $S_n = \sum_{j=1}^\infty
U_{n,j}^2$
to see
\[
{\cal L} ( S_n |\bfT_n=n\mu_n )
\Rightarrow{\cal L}\biggl (\int_0^1 \psi^2(s) W^2(s)\,\mathrm{d}s \biggr)
\]
which establishes (\ref{weakconvergenceb}) and completes the proof.
\end{enumerate}

\subsection{Unconditional limits}\label{sec2.2} \label{uncondweaklimitsec}
We now consider the random function
$
Y_n(t) = \psi(t) W_n(t)
$
and review some well known facts about the unconditional
limiting distributions of the processes $Y_n$; see
Shorack and Wellner \cite{shorackwellner}, for example.
If $\theta_n$ converges to $\theta$, then the unconditional laws
of $Y_n$ converge weakly in $D[0,1]$
to the law of a Gaussian process $Y$ with mean 0 and covariance
\[
\zeta_\theta(s,t) = \psi(s)\rho_\theta(s,t)\psi(t).
\]
The covariance
$\zeta_\theta$ is square integrable over the unit square;
it is convenient to suppress $\theta$ in the notation for what follows.
There is a sequence of bounded continuous orthonormal eigenfunctions
$ \chi_j(t),j=1,2,\ldots,$
with corresponding eigenvalues $\lambda_j$ such that
\[
\int_0^1 \zeta(s,t)\chi_j(t)\,\mathrm{d}t \equiv\lambda_j \chi_j(s).
\]
Then
%
\begin{equation}\label{sumofweightedchisquares}
\int_0^1 Y^2(t)\,\mathrm{d}t = \sum\lambda_j Z_j^2,
\end{equation}
where
\[
Z_j = \lambda_j^{-1/2} \int_0^1 Y(t) \chi_j(t)\,\mathrm{d}t.
\]
The $Z_j$ are independent standard normal.
Let $H_\infty(\cdot;\theta)$ denote the cumulative distribution
of (\ref{sumofweightedchisquares}). It is then standard that
\[
\lim_{n\to\infty}\sup_{-\infty<x<\infty} \sup_{\theta\in\Gamma}
|H_n(x;\theta)-H_\infty(x;\theta)| = 0.
\]
Our main result will therefore follow if we establish
(\ref{weakconvergenceb}).


Next, recall that $W_n$ converges weakly to the Gaussian process $W$
with covariance function
$\rho_\theta$. The map
\[
f\mapsto\biggl(\int_0^1 f(s) \psi(s) g_1(s)\,\mathrm{d}s,\ldots,\int_0^1 f(s) \psi
(s) g_K(s)\,\mathrm{d}s\biggr)
\]
is continuous from $D[0,1]$ to $\reals^K$ so that
\[
(U_{n,1},\ldots,U_{n,K}) \Rightarrow(U_{\infty,1},\ldots,U_{n,K}).
\]
This limit vector has a multivariate normal distribution with mean 0 and
covariance
%
\begin{equation}
\label{CovUlimit}
\operatorname{Cov}(U_{\infty,i},U_{\infty,j}) = \int_0^1\int_0^1 g_i(s)
g_j(t) \zeta_\theta(s,t)\,\mathrm{d}s \,\mathrm{d}t.
\end{equation}

It follows that
\[
\sum_{j=1}^K U_{n,j}^2 \Rightarrow\sum_{j=1}^K U_{\infty,j}^2.
\]
Since
\[
\int_0^1 \psi^2(s) W^2(s)\,\mathrm{d}s = \sum_{j=1}^\infty U_{\infty,j}^2
\]
almost surely we have, for any sequence $K_n$ tending to infinity, that
\[
\sum_{j=1}^{K_n} U_{\infty,j}^2\Rightarrow\int_0^1 \psi^2(s)
W^2(s)\,\mathrm{d}s.
\]

This completes the analysis of the unconditional limit behaviour of $S_n$.
The next subsection considers the conditional limit behaviour.

\subsection{Convergence of finite dimensional distributions --
conditional case}\label{sec2.3}

In the following, all distributional assertions are statements about the
conditional distribution of the objects involved given $\bfT_n = n\mu_n$
for a specific sequence $\theta_n$ converging to some $\theta\in
\Gamma$ and
$\mu_n = \nabla K (\theta_n)$.
We apply Lemma~\ref{Holst} as follows. We have
\[
U_{n,j} \equiv\int_0^1 \psi(t) W_n(t) g_j(t)\,\mathrm{d}t = n^{-1/2} \sum
_{i=1}^n\Phi_{jn}(X_i),
\]
where
\[
\Phi_{jn}(x) = \int_0^1  [1\{F(x;\theta_n)\le t\} - t ]
\psi(t) g_j(t)\,\mathrm{d}t.
\]
It follows from Lemma~\ref{Holst}
that
\[
{\cal L} \bigl((U_{n,1},\ldots,U_{n,K}) | \bfT_n=n\mu_n \bigr)
\Rightarrow
{\cal L} ((U_{\infty,1},\ldots,U_{\infty,K}) ).
\]
The vector $(U_{\infty,1},\ldots,U_{\infty,K})$ has a multivariate
normal distribution
with mean 0 and variance covariance matrix with entries as at~(\ref{CovUlimit}).
This is the same limit behaviour as in the unconditional case.
Thus,
\[
{\cal L} \Biggl( \sum_{j=1}^K U_{n,j}^2 \Big|\bfT_n=n\mu
_n \Biggr)
\Rightarrow\sum_{j=1}^K U_{\infty,j}^2.
\]
Again this is the same weak limit as in the previous section. Finally,
since convergence in distribution is metrizable there is a sequence $K_n$
tending to infinity so slowly that
\[
{\cal L} \Biggl( \sum_{j=1}^{K_n} U_{n,j}^2  \Big|\bfT_n=n\mu
_n \Biggr)
\Rightarrow\sum_{j=1}^\infty U_{\infty,j}^2.
\]

We need only show, therefore,
that for any sequence $K_n$ tending to infinity we have,
conditionally on $\bfT_n=n\mu_n$,
\[
\sum_{j=K_n+1}^\infty U_{n,j}^2 \to0
\]
in probability. It suffices to show that
%
\begin{equation} \label{tailexpectation}
{\rm E}\Biggl ( \sum_{j=K_n+1}^\infty U_{n,j}^2  \Big|\bfT
_n=n\mu_n
 \Biggr) \to0.
\end{equation}
We will prove this from the following statements.
First, we will show that for each fixed~$j$
%
\begin{equation} \label{Onej}
{\rm E} (U_{n,j}^2 | \bfT_n=n\mu_n ) \to{\rm E} \{
U_{\infty,j}^2 \}.
\end{equation}
This shows that for each fixed $K$ we have
%
\begin{equation} \label{Kj}
{\rm E} \Biggl( \sum_{j=1}^K U_{n,j}^2  \Big| \bfT_n=n\mu
_n \Biggr) \to
{\rm E} \Biggl(\sum_{j=1}^K U_{\infty,j}^2 \Biggr).
\end{equation}
Finally, we will show that
%
\begin{equation} \label{ALLj}
{\rm E} \Biggl( \sum_{j=1}^\infty U_{n,j}^2  \Big| \bfT_n=n\mu
_n \Biggr) \to
{\rm E} \Biggl(\sum_{j=1}^\infty U_{\infty,j}^2 \Biggr).
\end{equation}
Assertion~(\ref{tailexpectation}) is a straightforward consequence
of (\ref{Onej}) and (\ref{ALLj}).
It is now straightforward to apply Slutsky's theorem to complete
the proof of the main theorem.

Statements (\ref{Onej}) and (\ref{ALLj}) are proved in
Section~\ref{proofs}. The proofs relate
$ {\cal L} (U_{n,j}^2 | \bfT_n=n\mu_n ) $ to an integral involving
\[
{\rm E}\bigl (1(X_i \le x)1(X_k \le y)|\bfT_n=n\mu_n \bigr)
\]
and other similar Rao--Blackwell estimates.
They then use a conditional
Edgeworth expansion of Rao--Blackwell estimates which is of some
interest in
its own right. We describe these expansions in the next section.

\section{Conditional Edgeworth expansions}\label{sec3}

In this section, we compute the first term in an Edgeworth
expansion of the conditional expectation of a function of $X_1,\ldots,X_m$
given $\bfT_n$. We will focus on uniformity, extending the work of
Portnoy \cite{portnoy}. The calculations may be
interpreted as a computation of the difference, to order $1/n$, between a
Rao--Blackwell estimate of a parameter and the maximum likelihood estimate.

Our results use the Edgeworth expansion of the density of $\bfT_n$.
Assuming Condition~\ref{condD}, for $n\ge r$ the quantity
$\{\bfT_n-n\mu(\theta)\}/\sqrt{n}$ has a density $q_n(\cdot;\theta
)$. The
following lemma is essentially a uniform version of Theorem 19.2 in
Bhattacharya and Ranga~Rao \cite{bhatta};
see Holst \cite{holst}, Yuan and Clarke \cite{yuanclarke}.
It extends a lemma appearing in Lockhart and O'Reilly~\cite{LOreMoores}.
Let $\mathbf{u}$ denote a $k$ vector with entries $u_1,\ldots, u_k$.

\begin{lemma} \label{edgeworthlemma} Assume Condition~\textup{\ref{condD}}.
Then there are functions
\[
\psi_j(\bfu;\theta), \qquad  j=1,2, \ldots,
\]
and
\[
\psi_{jk}(\bfu;\theta), \qquad k=0,\ldots,j+2,
\]
such that
\begin{enumerate}[5.]
\item$\psi_{jk}$ is homogeneous of degree $k$ as a
function of $u_1,\ldots,u_k$.
That is
\[
\psi_{jk}(u_1,\ldots,u_k;\theta)= \sum_{i_1\cdots i_k}
a_{jk;i_1\cdots i_k}(\theta)
u_{i_1}\cdots u_{i_k}
\]
for some coefficients $ a_{jk;i_1\cdots i_k}(\theta)$ not depending on
$\bfu$.

\item If $j-k$ is odd, then $\psi_{jk}\equiv0$.

\item $\psi_j$ is a polynomial of degree $j+2$ as a function of $\bfu
$ given by
\[
\psi_j = \sum_{k=0}^{j+2} \psi_{jk}.
\]

\item The coefficients $a_{jk;i_1\cdots i_k}(\theta)$ in these
polynomials are smooth functions of $\theta$.

\item Fix an integer $s \ge0$ and a compact subset $\Gamma$ of
$\intTh$.
Let $\phi(\bfu,V)$ be the multivariate normal density with mean 0 and
covariance matrix $V$.
Then
\[
\epsilon_n\equiv
\sup_{\theta\in\Gamma}\sup_\bfu \Biggl| q_n(\bfu;\theta) -
\phi\{\bfu,V(\theta)\}\Biggl \{1+\sum_{j=1}^s
\frac{\psi_j(\bfu,\theta)}{n^{j/2}} \Biggr\} \Biggr|
= \mathrm{O}\bigl(n^{-(s+1)/2}\bigr).
\]
\end{enumerate}
\end{lemma}

We will use this lemma with $s=3$ to get an error rate on our 1 term expansion.
We need the following notation. Define
\[
B_m(x_1,\ldots,x_m) = \sum_{i=1}^m \{T(x_i) - \mu\}
\]
and let $\mathbf{B}_m$ denote the random vector
\[
B_m(X_1,\ldots,X_m)=\bfT_m-m\mu.
\]
Let $D = V^{-1}$ be the inverse of the variance covariance matrix $V$.
The lowest degree term in the polynomial $\psi_1$ has the form
\[
\psi_{1,1}(\bfu) = - \sum_1^k a_{1,1;\ell}(\theta) u_\ell
\]
where, from Bhattacharya and Ranga~Rao \cite{bhatta}, page 55,
we have
\begin{eqnarray*}
a_{1,1;\ell}& =& \sum_i   \kappa_{iii} D_{ii} D_{i\ell} / 2
+\sum_{i\neq j} \kappa_{iij}
 ( 2D_{ij} D_{i\ell} + D_{ii}D_{i\ell} )/2
\\
&&{} + \sum_{i < j < k} \kappa_{ijk}
 ( D_{ij} D_{k\ell} + D_{ik} D_{j\ell}+D_{jk}D_{i\ell} ).
\end{eqnarray*}
%

If $J(x_1,\ldots,x_m)$ is a real valued measurable function on $\Omega^m$;
we let $\mathbf{J} = J(X_1,\ldots,X_m)$.
Remember in the following that $\mu$ and $\theta$ are related
through ${\rm E}_\theta(\bfT_n)=n\mu$.

\begin{theorem} 
\label{edgeworthexpansion}
Fix an integer $m>0$.
Suppose ${\bar J}\ge0$ is a real valued measurable function on $\Omega^m$
such that
\[
{\rm E}_\theta \{{\bar J}(X_1,\ldots,X_m)
 \} < \infty
\]
for all $\theta\in\intTh$.
Then for each compact subset $\Gamma$ of $\intTh$ we have
%
\begin{equation} \label{remainderbound}
\limsup_{n\to\infty} n^2\sup_{\theta\in\Gamma} \sup_J  |
{\rm E}\{ \bfJ|\bfT_n=n\mu\}
-A(n,J,\theta) | < \infty,
\end{equation}
where
\[
A(n,J,\theta) \equiv{\rm E}_\theta (\mathbf{J} )+
\frac{R(J,\theta)}{n}
\]
and
%
\begin{eqnarray}
\label{R1def}
R(J,\theta) &\equiv&
\frac{mk}{2}{\rm E}_\theta (\mathbf{J} )-
\frac12{\rm E}_\theta \{\bfJ\mathbf{B}_m^\prime V^{-1}(\theta)
\mathbf{B}_m \}
- {\rm E}_\theta \{\bfJ\psi_{1,1}(\bfB_m) \}
\\
\label{R2def}
& =& \nabla^2 {\rm E}_\theta(\bfJ) + \psi_{1,1} \{\nabla{\rm
E}_\theta(\bfJ) \}.
\end{eqnarray}
The supremum over $J$ is over all measurable $J$ defined on $\Omega^m$ with
$|J|\le{\bar J}$ (almost everywhere).
Moreover,
%
\begin{equation} \label{term1bounded}
\sup_{\theta\in\Gamma} \sup_{J: |J| \le\bar{J}} R(J,\theta) =\mathrm{O}(1).
\end{equation}
\end{theorem}

In (\ref{R2def}), the symbols $\nabla$ and
$\nabla^2$ are as in Lemma~\ref{exponentialfamilyprops}.
It is part of the theorem that the quantities on the right
in (\ref{R1def}) and (\ref{R2def}) are equal.

\section{Examples}\label{sec4}

In this section, we consider the Gamma and von Mises models and show that
the theory of the previous sections applies. These two models were
considered in Lockhart, O'Reilly and Stephens \cite{LOS1,LOS2}
where Gibbs sampling was used to implement the conditional tests
discussed here
via Markov Chain Monte Carlo. In the case of the Gamma distribution, we
also illustrate the use of the expansion of the Rao--Blackwell estimate by
giving a~formula for an approximate Rao--Blackwell estimate of the
shape parameter.

\subsection{The Gamma distribution}\label{sec4.1}

Suppose $X_1,X_2,\ldots$ are i.i.d. with density
\[
f(x;\alpha,\beta) = \frac{1}{\beta\Gamma(\alpha)}
 \biggl(\frac{x}{\beta} \biggr)^{\alpha-1} \exp(-x/\beta) 1(x>0).
\]
We take $\theta_1=\alpha$, $\theta_2=1/\beta$ and
$\Theta=\{\theta:\theta_1>0,\theta_2>0\}$. We then
have
\[
T(x) = (\log(x), -x)\vadjust{\goodbreak}
\]
and
\[
\kappa(\phi_1,\phi_2) =
\log \biggl\{\frac{\Gamma(\theta_1+\phi_1)}{\Gamma(\theta_1)}
\frac%
{\theta_2^{\theta_1}}%
{(\theta_2+\phi_2)^{\theta_1+\phi_1}}%
 \biggr\}.
\]
The characteristic function of $\bfT$ is
\[
\Psi(\phi_1,\phi_2) =
\frac{\Gamma(\theta_1+\mathrm{i}\phi_1)}{\Gamma(\theta_1)}
\frac{(\theta_2+\mathrm{i}\phi_2)^{\theta_1+\mathrm{i}\phi_1}}{\theta_2^{\theta_1}}.
\]
Fix a compact set $\Gamma$ in the parameter space and let
\[
\epsilon= \inf\{\theta_1\dvt \exists\theta_2 \dvt (\theta_1,\theta
_2)\in
\Gamma\}.
\]
In Section~\ref{proofs}, we use properties of the Gamma function in
the complex plane
to show that for $r$ so large that $r\epsilon> 2$ and $r > 4$ we have
%
\begin{equation}\label{gammacondD}
\sup_{\theta\in\Gamma}\int
 |\Psi(\phi_1,\phi_2) )|^r \,\mathrm{d}\phi_1\,\mathrm{d}\phi_2 < \infty.
\end{equation}
This establishes Condition~\ref{condD} in this case.

For completeness, we record here the functions needed to apply Theorem
3 to this family. Let $\psi(\theta) = \mathrm{d}\log\Gamma(\theta)/\mathrm{d}\theta
$ denote the
digamma function and let
$\psi^\prime$ and $\psi^{\prime\prime}$ denote its first
and second derivatives.
Let $\delta= \theta_1\psi^{\prime}(\theta_1)-1$. Then we find
%
\begin{eqnarray*}
\mu_1 & =& \psi(\theta_1)-\log(\theta_2), \qquad
\mu_2 = -\theta_1/\theta_2,
\\
V_{11} &=& \psi^{\prime}(\theta_1)
, \qquad
V_{12} = V_{21}=1/\theta_2,
\\
V_{22} &=& \theta_1/\theta_2^2
, \qquad
D_{12}  =D_{21} = \theta_2/\delta,
\\
D_{11} &=& \theta_1/\delta
, \qquad
D_{22} = \theta_2^2\psi^{\prime}(\theta_1)/\delta,
\\ \kappa_{111} &=& \psi^{\prime\prime}(\theta_1) , \qquad  \kappa_{112}
= \kappa_{121} = \kappa_{211}=0, \\
\kappa_{222} &=& -2\theta_1/\theta_2^3 , \qquad  \kappa_{122} = \kappa
_{221} = \kappa_{221}=1/\theta_2^2,
\\
a_{11;1} &=& \frac{\theta_1^2\psi^{\prime\prime}(\theta_1)+2\theta
_2\psi^{\prime}(\theta_1) +2}{2\delta^2}
, \qquad
a_{11;2} =
\frac{\theta_1\psi^{\prime\prime}(\theta_1)+2\theta_2 \{
\psi^{\prime}(\theta_1) \}^2 +2\theta_2\psi^{\prime}(\theta
_1)}{2\delta^2}.
\end{eqnarray*}
These formulas may be used to give approximations in terms of the maximum
likelihood estimate $\hat\theta$ to order $1/n$
of the Rao--Blackwell estimate of a parameter.
As an example, we consider the approximation to the
Rao--Blackwell estimate of the shape parameter $\theta_1$. In this
case ${\rm E}_\theta(\bfJ) = \theta_1$ so the Hessian matrix in
$R(J,\theta)$
is 0 and the gradient is
simply $(1,0)^\prime$. Our approximation from (\ref{R2def}) is then
\[
\tilde\theta_1=\hat\theta_1-\frac{\hat\psi_{1,1}(1,0)}{n} =
\hat\theta_1+\frac{\hat\theta_1^2\psi^{\prime\prime}(\hat
\theta_1)+2\hat\theta_2\psi^{\prime\prime}(\hat\theta_1) +2}{%
2 n  \{\hat\theta_1\psi^{\prime}(\hat\theta_1)-1 \}^2}.
\]

\begin{remark*} I do not know if there is, for some value of $m$,
an unbiased estimate of $\theta_1$. That is, I do not know if\vadjust{\goodbreak} $\mathbf{J}$
exists in the calculation just given. It seems worth
noting that the expansion can be computed anyway since the terms
therein depend only on the function of
the parameters which is being estimated and the derivatives of that
function.\vspace*{-2pt}
\end{remark*}

\subsection{The von Mises distribution}\label{sec4.2} \label{vonmisessec}
\vspace*{-2pt}

Suppose $X_1,X_2,\ldots$ are i.i.d. with density
\[
f(x;\alpha,x_0) =
\frac{1}{2\uppi I_0(\alpha)}\exp\{\alpha\cos(x-x_0)\} 1( 0 < x <
2\uppi),
\]
where $I_0$ is the modified Bessel function of the first kind of order 0.
We take $\theta_1=\alpha\cos(x_0)$, $\theta_2=\alpha\sin(x_0)$ and
$\Theta=\reals^2$. We then
have
\[
T(x) = (\cos(x),\sin(x)).
\]
Here we find it easier to verify Condition~\ref{condD*}. For a
sample of size $m$ the
density of the sufficient statistics is known analytically in the case
$\theta_1=\theta_2=0$,
that is, when the distribution is uniform on the interval $(0,2\uppi)$. Write
$\bfT_m$ in polar coordinates as $(R\cos\delta, R\sin\delta)$ with
the angle $\delta$ in
$[0,2\uppi)$ and $R=\|\bfT_r\|$; then $R$ and $\delta$ are
independent. The distribution of
$\delta$ is uniform on $[0,2\uppi)$. From Stephens \cite{stephens62},
we find $R$ has the density
\[
f_m(u) = u\int_0^\infty J_0(ut)J_0^m(t) t \,\mathrm{d}t,
\]
where $J_0$ is the Bessel function of the first kind of order 0.
The function $J_0(t)$ is bounded and decays at infinity
like $t^{-1/2}$.
So
for all $m>4$
there is
a constant $C_m$ such that
\[
f_m(u) \le C_mu
\]
for all $u>0$. The density $f_m$ vanishes for negative $u$ and for $u>m$.
Change variables to see that for all $m\ge5$ the density of $\bfT_m$
is bounded by $C_m/(2\uppi)$. For $\theta=(\theta_1,\theta_2)$ not 0 the
likelihood ratio
of $\theta$ to $0$ is $\exp(\theta^\prime\bfT_m)/I_0^m(\|\theta\|)$.
Since the density of $\bfT_m$ for $\theta$ is the density for 0
multiplied by the likelihood ratio
Condition~\ref{condD*} holds with $r\equiv5$.\vspace*{-2pt}

\section{Discussion}\label{sec5}\label{discussion}
\vspace*{-2pt}

We conclude with a series of remarks.\vspace*{-2pt}

\begin{remark}
For a given goodness-of-fit test statistic we may compute $P$-values in
several ways.
The parametric bootstrap technique proceeds by estimating the unknown
parameters and then generating a large number of samples from the
hypothesized distribution using the estimated value of the parameters.
Except in location-scale models the resulting tests are approximate; that
is, the distribution of the $P$-value is not exactly uniform though it
becomes more so as the sample size increases.\vspace*{-2pt}
\end{remark}

An alternative technique is to compute a conditional $P$ value using
\[
P(S_n > s|\bfT_n)\vadjust{\goodbreak}
\]
evaluated at $s$ equal to the observed value of $S_n$.
This $P$-value must generally be evaluated by Monte Carlo methods. For some
distributions, such as the Inverse Gaussian, there is a
direct way to simulate samples from the conditional distribution of the
data given $\bfT_n$.
See
O'Reilly and Gracia-Medrano \cite{FJRandLGM}.
For other distributions, Markov Chain Monte Carlo may be used; see
Lockhart, O'Reilly and Stephens \cite{LOS1,LOS2}.

If the null hypothesis is true and the true value of $\theta$ is in
$\intTh$, then we have shown that the difference between these
two $P$-values converges almost surely to 0. In our experience,
these two $P$-values are usually extremely close together
suggesting the agreement extends to some higher order expansion; I do
not know how to show such a thing.

\begin{remark} Indeed this equivalence of $P$-values requires
only a large sample size and an estimate $\hat\theta$ not too close to
the boundary of $\Theta$. It is not at all necessary that the null
hypothesis be true. Of course if the null hypothesis is not true the
estimate $\hat\theta_n$ could converge to the boundary of the parameter
space and then our results permit the $P$-values to be different even
in large samples.
\end{remark}

\begin{remark} For fixed alternatives, our results imply
that the
difference in powers between the two tests tends to 0 except
when $\bfT_n/n$ does not have a limit in $\mu(\intTh)$.
The conclusions in Theorem~2 can be extended to contiguous sequences of
alternatives yielding conclusions that the two tests have identical
limiting powers along such sequences.
\end{remark}

\begin{remark}
The local central limit theorem for lattice distributions may be used
to prove the equivalent of Theorem~\ref{maintheorem} if $T(x)$ takes
values in a lattice and the data are discrete.
\end{remark}

\begin{remark} The result also extends to a variety
of other statistics such as
\[
\int_{0}^1  \biggl\{W_n(t) - \int_{0}^1\psi(u)W_n(u)\,\mathrm{d}u  \biggr\}^2
\psi^2(t)\,\mathrm{d}t
\]
or
\[
\int_0^1\int_0^1 K(s,t) W_n(s) W_n(t)\,\mathrm{d}s \,\mathrm{d}t
\]
or any other suitable quadratic form in the process $W_n$,
under regularity conditions on the weight functions $\psi$, the
kernel $K$, or the quadratic form.
\end{remark}

\begin{remark} One important case not covered
by our proof is the Anderson--Darling test which is of the
Cram\'{e}r--von~Mises type but with weight function
\[
\psi(s) = 1/\sqrt{s(1-s)}
\]
which is not square integrable. It may be possible
to verify our assertions (\ref{Onej}) and (\ref{ALLj})
by more careful
analysis of the conditional moments of $W_n$ near the ends of the unit
interval.
\end{remark}

\begin{remark} Our proofs show that the
Edgeworth expansion to order $2s$ given in  Lemma~\ref{edgeworthlemma} may be used to provide
an expansion of any Rao--Blackwell estimate about the maximum likelihood
estimate of ${\rm E}_\theta(\bfJ)$ in inverse powers of $n$ out to
terms of order $n^{-s}$ with a~remainder which is $\mathrm{O}(n^{-(s+1)})$
uniformly on compact subsets of $\intTh$.
We have not done the algebra for any $s>1$ but
we can state the following theorem.
\end{remark}
\begin{theorem}
Under the conditions of Theorem~\ref{edgeworthexpansion}, there are functions
$R_j(J,\theta)$ for $j=1,2, \ldots$ such that for any integer
$s \ge1$ we have
%
\begin{equation}
\limsup_{n\to\infty} n^{1+s}\sup_{\theta\in\Gamma} \sup_J  |
{\rm E}\{ \bfJ|\bfT_n=n\mu\}
-A_s(n,J,\theta) | < \infty,
\end{equation}
where
\[
A_s(n,J,\theta) \equiv{\rm E}_\theta (\mathbf{J} )+
\sum_{j=1}^s\frac{R_j(J,\theta)}{n^j}.
\]
The functions $R_j$ are computed using Taylor expansions as
in Theorem~\ref{edgeworthexpansion} and collecting
terms in inverse powers of $n$.
Each $R_i$ is bounded uniformly over $\theta\in\Gamma$ and $|J| \le
\bar{J}$.
\end{theorem}

Of course $R_1$ is just $R$ of Theorem~\ref{edgeworthexpansion} and
the point is that
the arguments in the proof of that theorem can be applied to all
remainder terms occurring here.

\begin{remark}\label{rem8} In Theorem~\ref{maintheorem},
the $X_i$ are real valued;
this is needed only for the weak convergence results.
In the von Mises case, for instance, it is useful to regard the
observation~$X_i$ not as an angle but as a unit vector~$X_i$ as was
suggested in the introduction.
This makes $\bfT_n = \sum X_i$. In many examples, the
$X_i$ can usefully be taken to be multivariate. Our
results
may be expected to extend to any statistic
admitting a sum of squares expansion like that of Cram\'er--von~Mises
statistics.
\end{remark}

\begin{remark}\label{rem9} The conditional tests described here have level
identically equal to $\alpha$. In the introduction, we noted that this is
a necessary condition for an unbiased level $\alpha$ test in models
with a
complete sufficient statistic. Though necessary, the condition is not
sufficient; we do not know how to check that a given conditional test is
unbiased, nor how to establish any optimal power properties for the tests
considered here.
\end{remark}

\section{Proofs}\label{sec6} \label{proofs}

\subsection{\texorpdfstring{Proof of Lemma~\protect\ref{Holst}}{Proof of Lemma 3}}\label{sec6.1}
The proof in Holst \cite{holst}
of his Corollary 3.6 extends directly to prove this lemma.
However, Holst's Corollary 3.6
assumes ``the general conditions'' of his Section~2. In particular,
we must verify the integrability hypothesis of his Proposition 2.1
which we now
describe in our notation. Let
\[
\Psi_{r,\theta}(\zeta_1,\zeta_2) =
{\rm E}_\theta\bigl\{\exp\bigl(\mathrm{i}\zeta_1^\prime S_r(\theta)+\mathrm{i}
\zeta_2^\prime\bfT_r\bigr)\bigr\}
\]
be the joint characteristic function of $S_r(\theta),\bfT_r$.
Holst requires that for
each $\zeta_1$ and each compact subset $\Gamma$ of $\intTh$
there is an $r>0$ such that
for all $\theta\in\Gamma$
%
\begin{equation} \label{Holstintegrability}
\int |\Psi_{r,\theta}(\zeta_1,\zeta_2) |\,\mathrm{d}\zeta_2 <
\infty.\vspace*{-2pt}
\end{equation}

\begin{lemma}
Condition~\textup{\ref{condD}} implies (\ref{Holstintegrability}). In fact, $r$ can be
chosen free of
$\zeta_1$.\vspace*{-2pt}
\end{lemma}

This is an easy consequence of the following lemma.\vspace*{-2pt}
\begin{lemma}
\label{HolstLemma}
Suppose $X\in\reals^n$ and $Y\in\reals^m$ have
joint distribution $F(\mathrm{d}x,\mathrm{d}y)$ and
joint characteristic function $\psi(u,v)$. Then
\begin{enumerate}[3.]
\item If $Y$ has density $f$ bounded by $M$ and $\psi$ is real valued
and nonnegative,
then
\[
\int_{\reals^m} \psi(u,v)\,\mathrm{d}v \le M(2\uppi)^m.
\]

\item If $Y$ has density $f$ bounded by $M,$ then
\[
\int_{\reals^m}  |\psi(u,v) |^2\,\mathrm{d}v \le M(2\uppi)^m.
\]

\item If
$
M \equiv\int_{\reals^m}|\psi(0,v)|\,\mathrm{d}v < \infty,
$
then $Y$ has a density $f$ such that for all $y$
\[
f(y) \le M/ (2\uppi)^m.\vspace*{-2pt}
\]
\end{enumerate}
\end{lemma}

\begin{pf}
Statement~3 is a well-known consequence of the Fourier inversion
formula. Statement~2 follows from Statement~1 by symmetrization: if
the pair $(X^*,Y^*)$ has the same joint distribution as $(X,Y)$ and is
independent of
$(X,Y)$ then the second statement is the first applied to $(X-X^*,Y-Y^*)$
noting that $Y-Y^*$ has a density also bounded by~$M$.

To prove Statement~1, we follow Feller \cite{feller2}, pages 480ff.
Let $\xi$ denote the standard normal density in $\reals^m$. Then for each
$a>0$ the function $a\xi(ax)$ is a density with characteristic function
$(2\uppi)^{m/2}\xi(u/a)$.
\begin{eqnarray*}
\int\exp\{ - \mathrm{i} \zeta^\prime v\} a^m \xi(av)\psi(u,v)\,\mathrm{d}v
& = &\int a^m \xi(av)\mathrm{e}^{\mathrm{i}u^\prime x} \exp\{\mathrm{i}v^\prime(y-\zeta)\}
F(\mathrm{d}x,\mathrm{d}y)\,\mathrm{d}v
\\[-2pt]
& = &(2\uppi)^{m/2}\int \mathrm{e}^{\mathrm{i}u^\prime x} \xi \{(y-\zeta)/a \}
F(\mathrm{d}x,\mathrm{d}y).
\end{eqnarray*}
At $\zeta=0$, we get
\begin{eqnarray*}
0 &\le&\int\psi(u,v) \exp(-a^2 v^\prime v/2)\,\mathrm{d}v
 \le (2\uppi)^{m}\int (1/a)^m \xi(y/a) F(\mathrm{d}x,\mathrm{d}y)
\\[-2pt]
& =& (2\uppi)^{m}\int (1/a)^m \xi(y/a) f(y)\,\mathrm{d}y
\le M(2\uppi)^m.
\end{eqnarray*}
Now let $a\to0 $ to get Statement~1.\vadjust{\goodbreak}
\end{pf}

\subsection{\texorpdfstring{Proof of Theorem~\protect\ref{edgeworthexpansion}}{Proof of Theorem 3}}\label{sec6.2}

We use the shorthands $\bfx$ for the vector $(x_1,\ldots,x_m)$ and
$\mathrm{d}\bfx$
for $\mu(\mathrm{d}x_1)\cdots\mu(\mathrm{d}x_m)$.
Let $f_m$ be the joint density of $X_1,\ldots,X_m$; we
suppress the dependence of this density on $\theta$. For $n \ge r$, we let
$q_n$ denote the density of
$ (\bfT_n-n\mu)/\sqrt{n}$
again suppressing the dependence on $\theta$. (Densities of
sufficient statistics are relative to Lebesgue measure while those of the
data are relative to products of the carrier measure $\mu$.) We adopt
the useful notation
\[
Q_m = B_m^\prime V^{-1} B_m, \qquad
Q_{mn}=Q_m/n
  \quad \mbox{and} \quad
q_n^*(x) = q_n(x)/\phi(0,V).
\]

It is elementary that
\[
{\rm E} \{\bfJ| \bfT_n=n\mu \} =
 \biggl(\frac{n}{n-m} \biggr)^{k/2}
\int J(\bfx) f_m(\bfx) \frac{q_{n-m}(A_m)}{q_n(0)}\,\mathrm{d}\bfx,
\]
where
\[
A_m = A_m(\bfx) =- \frac{\sum_{i=1}^m \{T(x_i) - \mu) \}
}{\sqrt{n}}= - \frac{B_m}{\sqrt{n}}.
\]

The quantity in~(\ref{remainderbound}) may be written as $| I_1+\cdots+I_8|$
where
$
I_i = \int J(\bfx)f_m(\bfx) \tau_i(\bfx)\,\mathrm{d}\bfx
$
for suitable functions $\tau_1,\ldots, \tau_{8}$. We will argue
below that each integral is $\mathrm{O}(n^{-2})$
uniformly in $\theta$ over compact subsets $\Gamma$ of $\intTh$.
The functions $\tau_i$ are given by
\begin{eqnarray*}
\tau_1(\bfu) & =& \nmk
\frac{q_{n-m}(A_m)- \phi(A_m,V) \{1+\sum_{j=1}^{4}
\fraca{\psi_j(A_m)}{(n-m)^{j/2}} \}}{q_n(0)},
\\
\tau_2(\bfu) & =&  \biggl\{\nmk-\nmkalt \biggr\} \frac{\phi
(A_m,V) \{1+\sum_{j=1}^{4} \fraca{\psi
_j(A_m)}{(n-m)^{j/2}} \}}{q_n(0)},
\\
\tau_3(\bfu) & =&\nmkalt \biggl[\frac{1}{q_n^*(0)} - \biggl\{1-\frac
{\psi_{2,0}}{n} \biggr\}\biggr ]
\mathrm{e}^{-Q_{mn}/2} \Biggl\{1+\sum_{j=1}^{4} \frac{\psi
_j(A_m)}{(n-m)^{j/2}} \Biggr\},
\\
\tau_4(\bfu) & =&\nmkalt \biggl(1-\frac{\psi_{2,0}}{n}
\biggr)\mathrm{e}^{-Q_{mn}/2} \biggl\{\sum_{j+\ell\ge4} \frac{\psi_{j\ell
}(A_m)}{(n-m)^{j/2}} \biggr\}
\\
& =& \nmkalt  \biggl(1-\frac{\psi_{2,0}}{n} \biggr)\mathrm{e}^{-Q_{mn}/2}
\biggl\{\sum_{j+\ell\ge4} \frac{(-1)^\ell\psi_{j\ell}(B_m)}{n^{\ell
/2}(n-m)^{j/2}} \biggr\},
\\
\tau_5(\bfu) & =& \nmkalt  \biggl(1-\frac{\psi_{2,0}}{n}
\biggr)\mathrm{e}^{-Q_{mn}/2}\psi_{1,1}(B_m)
 \biggl\{\frac{1}{n} - \frac{1}{\sqrt{n(n-m)}}  \biggr\},
\\
\tau_6(\bfu) & =& \nmkalt  \biggl(1-\frac{\psi_{2,0}}{n}
\biggr)\mathrm{e}^{-Q_{mn}/2}\psi_{2,0}\biggl \{\frac{1}{n-m} -\frac{1}{n}
\biggr\},
\\
\tau_7(\bfu) & =& \nmkalt
 \biggl(1-\frac{\psi_{2,0}}{n} \biggr)
\biggl (\mathrm{e}^{-Q_{mn}/2}-1+\frac{Q_m}{2n} \biggr)
 \biggl\{1+\frac{\psi_{2,0} - \psi_{1,1}(B_m)}{n} \biggr\},
\\
\tau_8(\bfu) & =& \nmkalt
 \biggl(1-\frac{\psi_{2,0}}{n} \biggr)
 \biggl(1-\frac{Q_m}{2n} \biggr)
 \biggl\{1+\frac{\psi_{2,0}-\psi_{1,1}(B_m)}{n} \biggr\} - \tau
_9(\bfu),
\end{eqnarray*}
where
\[
\tau_9(\bfu) = 1+\frac{ mk/2- Q_m/2 -\psi_{1,1}(B_m)}{n}.
\]
Theorem~\ref{maintheorem} will follow if we show for $i=1,\ldots, 8$ that
\[
\sup_{|J| \le{\bar J}}\sup_{\theta\in\Gamma} |I_i| = \mathrm{O}(n^{-2}).
\]
These 8 assertions may be established using several bounds. We do not
give complete details since the arguments are routine but we illustrate
some of the details.
For instance, it is elementary that
\[
\nmk \le(m+1)^{k/2}
 \quad \mbox{and} \quad
\nmk-\nmkalt=\mathrm{O}(n^{-2}).
\]
Continuity and compactness imply
\[
\sup_{\theta\in\Gamma}
\sup_x
 \Biggl|
\phi(x,V) \Biggl\{1+\sum_{j=1}^s\frac{\psi_j(x)}{n^{j/2}} \Biggr\}
 \Biggr| < \infty
\]
and
\[
\inf_{\theta\in\Gamma} \phi(0,V) > 0.
\]
Lemma~\ref{exponentialfamilyprops} guarantees that
\[
\liminf_{n\to\infty}\inf_{\theta\in\Gamma}q_n(0) > 0
\]
and so with $\epsilon_n$ as in Lemma~\ref{edgeworthlemma} we have
\[
|I_1| \le (m+1)^{k/2} \epsilon_n
\sup_{\theta\in\Gamma}{\rm E}_\theta ({\bar\mathbf{J}} )\big/
\inf_{\theta\in\Gamma}q_n(0).
\]
For $I_2$, $I_5$ and $I_6$ use the elementary facts that
\[
\frac{1}{n-m} - \frac{1}{n} =\mathrm{O}(n^{-2})
 \quad \mbox{and} \quad
\frac{1}{\sqrt{n(n-m)}} - \frac{1}{n} =\mathrm{O}(n^{-2}).
\]
Integral $I_3$ is bounded using Lemma~\ref{exponentialfamilyprops} again.
Integral $I_4$ uses the powers of $n$ in the displayed sum.
For $I_7$ use the inequalities
$ 0 < \mathrm{e}^{-x} - 1 +x < x^2/2 $
to see that
\[
0 <  \biggl(\mathrm{e}^{-Q_{mn}/2}-1+\frac{Q_m}{2n} \biggr) <
\frac{Q_m^2}{4n^2}.
\]
These bounds apply to the integrands; they are used to bound the
integrals based on the following observation.
The condition that $\bar{\mathbf{J}}$ have finite expectation for
all $\theta$ in $\intTh$ means that
$
{\bar J}(\mathbf{x}) f_m(\bfx)/{\rm E}_\theta({\bar\mathbf{J}})
$
defines another exponential family with natural parameter space
including $\intTh$. This permits differentiation under the integral sign
with respect to $\theta$ as many times as desired. It is then easily
established that for all $\alpha>0$
\[
\sup_{\theta\in\Gamma} {\rm E}_\theta (\|\bfT_r\|^\alpha
{\bfJ} ) < \infty.
\]
This permits all the bounds derived above to be integrated against
$J(\bfx)f_m(\bfx)$
to establish the desired conclusion.

Differentiation under the integral sign permits us to show
for any $J$ with $|J|\le{\bar J}$ the following two identities:
\begin{eqnarray*}
\nabla{\rm E}_\theta(\bfJ) &=& \operatorname{Cov}_\theta(\bfJ,\bfT_m),
\\
\nabla^2 {\rm E}_\theta(\bfJ) &=& \operatorname{Cov}_\theta (\bfJ,\bfB
_m\bfB_m^\prime )
\\
& =& {\rm E}_\theta (\bfJ\bfB_m \bfB_m^\prime )-{\rm
E}_\theta(\bfJ)V.
\end{eqnarray*}
From these two identities, we deduce
\begin{eqnarray*}
{\rm E}_\theta (\bfJ\bfB_m^\prime V^{-1} \bfB_m ) & =&
\operatorname{trace} \{{\rm E}_\theta (\bfJ\bfB_m \bfB_m^\prime
 ) V^{-1} \}
\\
& = &\operatorname{trace} \{\nabla^2 {\rm E}_\theta(\bfJ)V^{-1} \}
+ {\rm E}_\theta(\bfJ)\operatorname{trace }  (V^{-1} V ).
\end{eqnarray*}
This and the observation that $\psi_{1,1}$ is a linear function
establish the equivalence of the two forms of
$R(J,\theta)$ in (\ref{R1def}) and (\ref{R2def}).

\subsection{\texorpdfstring{Proof of assertions (\protect\ref{Onej}) and (\protect\ref{ALLj})}{Proof of assertions (9) and (10)}}\label{sec6.3}
\label{mainproof}

We must prove
\[
\label{Uimean}
{\rm E} [ U_{nj}^2|\bfT_n=n\mu_n ] \to\int_0^1\int_0^1
\psi(s)\psi(t) g_j(s) g_j(t) \rho_\theta(s,t)\,\mathrm{d}s\,\mathrm{d}t
\]
and
\[
\label{Umean}
{\rm E} [ S_{n}|\bfT_n=n\mu_n ] \to\int_0^1\psi
^2(s)\rho_\theta(s,s)\,\mathrm{d}s.
\]
To this end, define
\[
\tilde F(u|\mu) = {\rm E} \{1(X_1 \le x)| \bfT_n=n\mu \},
\]
where $u$ is related to $x$ by
$u=F(x,\theta)$.
Then ${\tilde F}(u|\bfT_n/n)$ is the Rao--Blackwell estimate of
$F(x,\theta)$.
Also define
$u_i=F(x_i,\theta)$ for $i=1,2$ and
\[
\tilde F(u_1,u_2|\mu) = {\rm E} \{1(X_1 \le x_1,X_2 \le x_2)|
\bfT_n=n\mu \}.
\]
Then ${\tilde F}_2(u_1,u_2|\bfT_n/n)$ is the Rao--Blackwell estimate of
$F(x_1,\theta)F(x_2,\theta)$
(the unconditional joint cumulative distribution function of $X_1$ and
$X_2$).

Define
\[
\rho_n(u_1,u_2|\mu) = {\rm E} \{W_n(u_1)W_n(u_2)| \bfT_n=n\mu
 \}.
\]
We then have
\[
{\rm E} \biggl[ \biggl \{\int_0^1 Y_n(t)g_j(t)\,\mathrm{d}t  \biggr\}
^2 \big|\bfT_n=n\mu_n  \biggr]
=
\int_0^1\int_0^1 \psi(s)\psi(t) g_j(s) g_j(t) \rho_n(s,t|\mu)\,\mathrm{d}s\,\mathrm{d}t.
\]
Direct calculation shows that
%
\begin{eqnarray}
\rho_n(u_1,u_2|\mu) &=&
{\tilde F}(\min(u_1,u_2)|\mu)
-{\tilde F}(u_1|\mu) u_2
-{\tilde F}(u_2|\mu) u_1
+u_1u_2 \nonumber
\nonumber\\
&&{}  +
(n-1) \{
{\tilde F}_2(u_1,u_2|\mu)
-{\tilde F}(u_1|\mu) u_2
-{\tilde F}(u_2|\mu) u_1
+u_1u_2 \}
\nonumber\\\label{CovTerm1}
& =&
{\tilde F}(\min(u_1,u_2)|\mu) -
{\tilde F}(u_1|\mu){\tilde F}(u_2|\mu)
\\
&&{}  +
(n-1)  \{
{\tilde F}_2(u_1,u_2|\mu) - {\tilde F}(u_1|\mu) {\tilde F}(u_2|\mu)
 \} 
\nonumber
\\\label{CovRemainder}
&&{}  + n
 \{ {\tilde F}(u_1|\mu)- u_1 \}
 \{ {\tilde F}(u_2|\mu)- u_2 \}.
\end{eqnarray}

We will establish (\ref{Onej}) by proving
%
\begin{equation}
\label{rholimit}
\rho_n(u_1,u_2|\mu) \to\rho_\theta(u_1,u_2)
\end{equation}
uniformly in $u_1$ and $u_2$.
We apply
Theorem~\ref{edgeworthexpansion}. Take ${\bar J}\equiv1$,
$J_1(X_1,X_2) = 1(X_1 \le x_1)$, $J_2(X_1) = 1(X_1 \le x_2)$ and
$J_3(X_1,X_2) = 1(X_1 \le x_1,X_2 \le x_2)$.
(The odd looking indexes in $J_2$ are deliberate.
The algebra involved in simplifying the remainder terms is easier if we
take $m=2$ for $J_3$ and $m=1$ for $J_1$ and $J_2$.)
We find from (\ref{term1bounded}) applied to $J_1$ and~$J_2$
that the term (\ref{CovRemainder}) converges to 0 uniformly in $u_1$ and
$u_2$. Applying (\ref{term1bounded}) to $J_1$ shows that the
term (\ref{CovTerm1}) converges, uniformly in $u_1$ and $u_2$, to
\[
\min(u_1,u_2)-u_1u_2.
\]

Finally from (\ref{remainderbound}), we find that
\[
(n-1)  \{
{\tilde F}_2(u_1,u_2|\mu) - {\tilde F}(u_1|\mu) {\tilde F}(u_2|\mu)
 \}
\]
converges to
\[
A(n,J_3,\theta)-A(n,J_1,\theta)A(n,J_2,\theta)
\]
uniformly in $u_1,u_2$.
Adopt the temporary notation
$
R_i = R(J_i,\theta)
$
and
$
A_i = A(n,J_i,\theta)
$
for $i=1,2,3$. Then
%
\begin{equation}
\label{approxlimit}
n (A_3-A_1A_2 )= R_3-R_1{\rm E}_\theta(\bfJ_2) - R_2 {\rm
E}_\theta(\bfJ_1) + R_1R_2/n.
\end{equation}
From (\ref{term1bounded}), we see that $R_1R_2/n$ converges to 0
uniformly in
$u_1$, $u_2$, $x_1$ and $x_2$.

Computing we get
\begin{eqnarray*}
R_3 & =& k{\rm E}_\theta(\bfJ_1){\rm E}_\theta(\bfJ_2)
-\frac12{\rm E}_\theta(\bfJ_3\mathbf{B}_2V^{-1}\mathbf{B}_2) + {\rm
E}_\theta(\bfJ_3\psi_{1,1}(\bfB_2)),
\\[-2pt]
R_1 {\rm E}_\theta(\bfJ_2)& =& \frac{k}{2}{\rm E}_\theta(\bfJ
_1){\rm E}_\theta(\bfJ_2)
-\frac12{\rm E}_\theta(\bfJ_1\mathbf{B}_1V^{-1}\mathbf{B}_1) + {\rm
E}_\theta(\bfJ_1\psi_{1,1}(\bfB_1)),
\\[-2pt]
R_2 {\rm E}_\theta(\bfJ_1)& =& \frac{k}{2}{\rm E}_\theta(\bfJ
_1){\rm E}_\theta(\bfJ_2)
-\frac12{\rm E}_\theta(\bfJ_2\mathbf{B}_1V^{-1}\mathbf{B}_1) + {\rm
E}_\theta(\bfJ_2\psi_{1,1}(\bfB_1)).
\end{eqnarray*}
Since $\bfB_2$ is a sum of two independent terms we expand the
quadratic form in $R_3$ to see
\begin{eqnarray*}
{\rm E}_\theta(\bfJ_3\mathbf{B}_2V^{-1}\mathbf{B}_2) &=&
{\rm E}_\theta(\bfJ_1\mathbf{B}_1V^{-1}\mathbf{B}_1) {\rm E}_\theta(\bfJ
_2) \\ &&{}+
{\rm E}_\theta(\bfJ_2\mathbf{B}_1V^{-1}\mathbf{B}_1) {\rm E}_\theta(\bfJ
_1) +
2{\rm E}_\theta(\bfJ_1\bfB_1^\prime)V^{-1} {\rm E}_\theta(\bfB
_1\bfJ_2).
\end{eqnarray*}
We may also use the linearity of $\psi_{1,1}$ and the independence of
$X_1$ and $X_2$ to see that
\[
{\rm E}_\theta(\bfJ_3\psi_{1,1}(\bfB_2)) = {\rm E}_\theta(\bfJ
_1\psi_{1,1}(\bfB_1)){\rm E}_\theta(\bfJ_2) +
{\rm E}_\theta(\bfJ_2\psi_{1,1}(\bfB_1)){\rm E}_\theta(\bfJ_1).
\]
Thus,
$
R_3-R_1{\rm E}_\theta(\bfJ_2) - R_2 {\rm E}_\theta(\bfJ_1) $
simplifies to
$
-{\rm E}_\theta(\bfJ_1\bfB_1^\prime)V^{-1} {\rm E}_\theta(\bfB
_1\bfJ_2).
$
Since $V$ is the Fisher information matrix in this problem, we have
established (\ref{Onej}).
To check (\ref{ALLj}), we make a very similar calculation.

\subsection{\texorpdfstring{Verification of Condition~\protect\ref{condD} for the Gamma family}{Verification of Condition D for the Gamma family}}\label{sec6.4}

Here, we establish (\ref{gammacondD}).
Change variables via $u=\phi_2/\theta_2$ to show the integral
in (\ref{gammacondD})
is proportional to $\theta_2^r$; thus we take $\theta_2=1$ without loss.
The integral becomes:
\[
\sup_{\theta\in\Gamma}\int\int \biggl|
\frac{\Gamma(\theta_1+\mathrm{i}\phi_1)}{\Gamma(\theta_1)} \biggr|^r
\frac{1}{(1+\phi_2^2)^{r\theta_1/2}}
\exp\{r\phi_1\tan^{-1}\phi_2\}\,\mathrm{d}\phi_1\,\mathrm{d}\phi_2.
\]
The substitution $\phi_2=\tan(u)$ reduces the integral to
\[
\sup_{\theta\in\Gamma}\int_{-\infty}^\infty\int_{-\uppi/2}^{\uppi/2}
 \biggl| \frac{\Gamma(\theta_1+\mathrm{i}\phi_1)}{\Gamma(\theta_1)} \biggr|^r
\cos^{r\theta_1-2}(u)\exp(r\phi_1u)\,\mathrm{d}u
\,\mathrm{d}\phi_1.
\]
We integrate separately over 4 ranges:
$R_1=\{-M \le\phi_1 \le M \}$,
$R_2=\{|\phi_1|>M, \phi_1 u < 0\}$,
$R_3=\{\phi_1>M, u>0\}$ and
$R_4=\{\phi_1< -M, u < 0\}$.
Since
$|\Gamma(\theta_1+\mathrm{i}\phi_1)| =|\Gamma(\theta_1-\mathrm{i}\phi_1)|$
the integrals $R_3$ and $R_4$ are equal.
Over $R_1$ we use the inequality
\[
|\Gamma(\theta_1+\mathrm{i}\phi_1)/\Gamma(\theta_1)| \le1
\]
(because the quantity inside the modulus signs is the characteristic function
of $\log(X_1)$)
to get the bound, for $\theta_1\ge\epsilon$ with $r\epsilon>2$
\begin{eqnarray*}
\int_{R_1}
 \biggl| \frac{\Gamma(\theta_1+\mathrm{i}\phi_1)}{\Gamma(\theta_1)} \biggr|^r
\cos^{r\theta_1-2}(u)
\mathrm{e}^{r\phi_1u}\,\mathrm{d}u
\,\mathrm{d}\phi_1 &\le& M \exp\{Mr\uppi/2\}\int_{-\uppi/2}^{\uppi/2}
\cos^{r\epsilon-2}(u)\,\mathrm{d}u
\\[-2pt]
& \le&\uppi M \exp\{Mr\uppi/2\}.
\end{eqnarray*}

Over $R_2$ the term $\exp(r\phi_1 u)$ is bounded by 1.
Thus,
\[
\int_{R_2}
 \biggl| \frac{\Gamma(\theta_1+\mathrm{i}\phi_1)}{\Gamma(\theta_1)} \biggr|^r
\cos^{r\theta_1-2}(u)
\mathrm{e}^{r\phi_1u}\,\mathrm{d}u
\,\mathrm{d}\phi_1 \le\uppi\int_{-\infty}^\infty
 \biggl| \frac{\Gamma(\theta_1+\mathrm{i}\phi_1)}{\Gamma(\theta_1)}
\biggr|^r \,\mathrm{d}\phi_1.
\]
The integral is bounded by the supremum of the density of $\log(X_1)$
over the real line and the compact parameter set $\Gamma$.

Finally, we consider the integral over $R_3$.
From Section 6.1.45 of Abramowitz and Stegun \cite{abram},
we find there is a constant $C$ such that
\[
 \biggl|\frac{\Gamma(\theta_1+\mathrm{i}\phi_1)}{\Gamma(\theta_1)} \biggr|
\le C \mathrm{e}^{-\uppi\phi_1/2} \phi_1^{\theta_1-1/2}.
\]
For $\theta_1 < 1/2 $, $\phi_1>M\ge1$ and $r\epsilon>2$ we then get
\begin{eqnarray*}
\int_{R_3}
 \biggl| \frac{\Gamma(\theta_1+\mathrm{i}\phi_1)}{\Gamma(\theta_1)} \biggr|^r
\cos^{r\theta_1-2}(u)
\mathrm{e}^{r\phi_1u}\,\mathrm{d}u
\,\mathrm{d}\phi_1
&\le& C
\int_0^{\uppi/2}
\int_{M}^\infty
\mathrm{e}^{-\phi_1 r(\uppi/2-u)} \cos^{r\theta_1-2}(u)\,\mathrm{d}\phi_1 \,\mathrm{d}u
\\[2pt] & \le& C\int_0^{\uppi/2} \frac{\cos^{r\theta_1-2}(u)}{r(\uppi/2-u)}\,\mathrm{d}u
\\[2pt]
&=& C\int_0^{\uppi/2} \frac{\sin^{r\theta_1-2}(u)}{ru}\,\mathrm{d}u
\\[2pt]
& \le& C\int_0^{\uppi/2} \frac{\sin^{r\epsilon-2}(u)}{ru}\,\mathrm{d}u
\\[2pt] & \le&\frac{C}{r} \int_0^{\uppi/2} u^{r\epsilon-3}\,\mathrm{d}u < \infty.
\end{eqnarray*}
For $\theta_1 \ge1/2$ we get
\begin{eqnarray*}
&&\int_{R_3}
 \biggl| \frac{\Gamma(\theta_1+\mathrm{i}\phi_1)}{\Gamma(\theta_1)} \biggr|^r
 \cos^{r\theta_1-2}(u)\exp(r\phi_1u)\,\mathrm{d}u
\,\mathrm{d}\phi_1
\\[2pt]
&& \quad  \le C
\int_0^{\uppi/2}
\int_{0}^\infty
\mathrm{e}^{-\phi_1 r(\uppi/2-u)} \phi_1^{r(\theta_1-1/2)} \cos^{r\theta
_1-2}(u)\,\mathrm{d}\phi_1 \,\mathrm{d}u
\\[2pt] && \quad
\le C\frac{\Gamma(1+r(\theta_1-1/2))}{r^{1+r(\theta_1-1/2)}}
\int_0^{\uppi/2} \frac{\sin^{r\theta_1-2}(u)}{u^{r(\theta
_1-1/2)+1}}\,\mathrm{d}u
\\[2pt] && \quad  \le
C\frac{\Gamma(1+r(\theta_1-1/2))}{r^{1+r(\theta_1-1/2)}}
\int_0^{\uppi/2} u^{r/2-3}\,\mathrm{d}u.
\end{eqnarray*}
For $r \ge5$ the right hand side is uniformly bounded over
$\Gamma\cap\{\theta_1 \ge1/2\}$.

\section*{Acknowledgements}
The author is grateful to anonymous reviewers whose advice lead to what
he hopes is a
much clearer paper.
He also thanks Federico O'Reilly and Michael
Stephens for many useful conversations on the topics discussed here.
The author acknowledges grant support from the
Natural Sciences and Engineering Research Council of Canada.

%

\printhistory

\end{document}